\documentclass[11pt]{article}
\usepackage{}
\usepackage{amssymb}
\usepackage{epsfig}

\def\be{\begin{equation}}
\def\ee{\end{equation}}
\def\bea{\begin{eqnarray}}
\def\eea{\end{eqnarray}}
\def\bes{\begin{eqnarray*}}
\def\ees{\end{eqnarray*}}

\def\nn{\nonumber}
\def\<{\langle}
\def\>{\rangle}
\def\lb{\label}

\def\R{{\bf R}}

\def\Z{{\bf Z}}
\def\N{{\bf N}}

\def\Q{{\bf Q}}

\def\Sg{{\Sigma}}

\def\hb{\vrule height0.18cm width0.14cm $\,$}

\title{Stable P-symmetric closed characteristics on partially symmetric compact convex hypersurfaces}
\author{Hui Liu$^{1}$, \thanks{Partially supported by NSFC (No.11401555), China Postdoctoral Science Foundation
  No.2014T70589, CUSF (No. WK0010000037). E-mail:huiliu@ustc.edu.cn.}
\qquad  Duanzhi Zhang$^{2}$, \thanks{Partially supported by NSFC (No.11422103, 11171341, 11271200) and LPMC of Nankai
University. E-mail: zhangdz@nankai.edu.cn}\\ \\
$^{1}$ Wu Wen-Tsun Key Laboratory of Mathematics, USTC, Chinese Academy of Sciences,\\
School of Mathematical Sciences, University of Science and Technology of China,
\\Hefei, Anhui 230026, P. R. China \\
$^{2}$ School of Mathematical Sciences and LPMC, Nankai University,\\ Tianjin 300071, P. R. China \\\\
Dedicate to Professor Rou-Huai Wang's 90th birth anniversity}
\date{}

\begin{document}

\maketitle

\begin{abstract}
{\it  In this paper, let $n\geq2$ be an integer,
$P=diag(-I_{n-\kappa},I_\kappa,-I_{n-\kappa},I_\kappa)$ for some
integer $\kappa\in[0, n-1)$, and $\Sigma \subset {\bf R}^{2n}$ be a
partially symmetric compact convex hypersurface, i.e., $x\in \Sigma$ implies $Px\in\Sigma$.
We prove that if $\Sigma$ is $(r,R)$-pinched with
$\frac{R}{r}<\sqrt{\frac{5}{3}}$, then $\Sigma$ carries at least two geometrically distinct P-symmetric closed
characteristics which possess at least $2n-4\kappa$ 	
Floquet multipliers on the unit circle of the complex plane.}
\end{abstract}



\renewcommand{\theequation}{\thesection.\arabic{equation}}
\renewcommand{\thefigure}{\thesection.\arabic{figure}}

\setcounter{equation}{0}
\section{Introduction and main results}
In this paper, we consider the stability of P-symmetric
closed characteristics on partially symmetric hypersurfaces in ${\bf R}^{2n}$.
Let $\Sigma$ be a $C^3$ compact hypersurface in ${\bf R}^{2n}$,
bounding a strictly convex compact set $U$ with non-empty interior, where $n\geq2$.
We denote the set of all such hypersurfaces by $\mathcal{H}(2n)$.
Without loss of generality, we suppose $U$ contains the origin. We
consider closed characteristics $(\tau,y)$ on $\Sigma$, which are
solutions of the following problem \be
\left\{\matrix{\dot{y}(t)=JN_\Sigma(y(t)),~y(t)\in \Sg,~\forall~t\in
{\bf R}, \cr
               y(\tau)=y(0), ~~~~~~~~~~~~~~~~~~~~~~~~~~~~~~~~~\cr }\right. \lb{1.1}\ee
where $J=\left(
              \begin{array}{cc}
                0 & -I_{n} \\
                I_{n} & 0 \\
              \end{array}
            \right)
$, $I_{n}$ is the identity matrix in ${\bf R}^n$ and
$\mathit{N}_\Sigma(y)$ is the outward normal unit vector  of
$\Sigma$ at $y$ normalized by the condition
$\mathit{N}_\Sigma(y)\cdot y=1$. Here $a\cdot b$ denotes the
standard inner product of $a, b \in {\bf R}^{2n}$. A closed
characteristic $(\tau, y)$ is {\it prime} if $\tau$ is the minimal
period of $y$. Two closed characteristics $(\tau,x)$ and
$(\sigma,y)$ are {\it geometrically distinct}, if $x({\bf
R})\not=y({\bf R})$. We denote by $\mathcal{J}(\Sigma)$ the set of all closed characteristics
$(\tau,y)$ on $\Sigma$ with $\tau$ being the minimal period of $y$.
For any $s_i, t_i\in{\bf R}^{k_i}$ with $i=1,2$, we denote by $(s_1,
t_1)\diamond(s_2, t_2)=(s_1, s_2, t_1, t_2)$. Fixing an integer
$\kappa$ with $0\leq \kappa < n-1$, let
$P=diag(-I_{n-\kappa},I_\kappa,-I_{n-\kappa},I_\kappa)$ and
$\mathcal {H}_\kappa(2n)=\{\Sigma\in\mathcal {H}(2n)\mid x\in\Sigma
~implies~ Px\in\Sigma\}$. For $\Sigma\in\mathcal {H}_\kappa(2n)$,
let $\Sigma(\kappa)=\{z\in {\bf R}^{2\kappa}\mid 0\diamond
z\in\Sigma\}$, where 0 is the origin in ${\bf R}^{2n-2\kappa}$. As
in \cite{DoL1}, A closed characteristic $(\tau, y)$ on $\Sigma \in
\mathcal {H}_\kappa(2n)$ is {\it P-asymmetric} if $y({\bf R})\cap
Py({\bf R})=\emptyset$, it is {\it P-symmetric} if $y({\bf
R})=Py({\bf R})$ with $y=y_1\diamond y_2$ and $y_1\neq 0$, or it is
{\it P-fixed} if $y({\bf
R})=Py({\bf R})$ and $y=0\diamond y_2$, where
$y_1\in{\bf R}^{2(n-\kappa)}$, $y_2\in{\bf R}^{2\kappa}$. We call a
closed characteristic $(\tau, y)$ is {\it P-invariant} if $y({\bf
R})=Py({\bf R})$. Then a P-invariant
closed characteristic is P-symmetric or
P-fixed.

Let $j: {\bf R}^{2n}\to {\bf R}$ be the gauge function of $\Sigma$,
i.e., $j(\lambda x)=\lambda$ for $x\in\Sigma$ and $\lambda\geq 0$,
then $j\in C^{3}({\bf R}^{2n}\setminus\{0\},{\bf R})\cap C^{0}({\bf
R}^{2n},{\bf R})$ and $\Sigma=j^{-1}(1)$. Fix a constant $\alpha
\in(1, +\infty)$ and define the Hamiltonian $H_{\alpha}: {\bf
R}^{2n}\to[0,+\infty)$ by\[ H_{\alpha}(x) :=j(x)^{\alpha}\] Then
$H_{\alpha}\in C^{3}({\bf R}^{2n}\setminus\{0\},{\bf R})\cap
C^{0}({\bf R}^{2n},{\bf R})$ is convex and
$\Sigma=H_{\alpha}^{-1}(1)$. It is well known that the problem
$(1.1)$ is equivalent to the following given energy problem of the
Hamiltonian system
$$\left\{\begin{array}{ll}
\dot{y}(t)=JH_{\alpha}^\prime(y(t)), H_{\alpha}(y(t))=1,~\forall~t\in {\bf R},\\
y(\tau)=y(0).
\end{array}\right. \eqno(1.2)$$Denote by $\mathcal{J}(\Sigma,\alpha)$ the set of all solutions
$(\tau,y)$ of the problem $(1.2)$, where $\tau$ is the minimal
period of $y$. Note that elements
in $\mathcal{J}(\Sigma)$ and $\mathcal{J}(\Sigma,\alpha)$ are in one
to one correspondence with each other. Let $(\tau,y)
\in\mathcal{J}(\Sigma,\alpha)$. We call the fundamental solution
$\gamma_y: [0,\tau]\to Sp(2n)$ with $\gamma_y(0)=I_{2n}$ of the
linearized Hamiltonian system\[
\dot{z}(t)=JH_{\alpha}^{\prime\prime}(y(t))z(t),~\forall~t \in {\bf
R}.\] the {\it associated symplectic path} of $(\tau,y)$. The
eigenvalue of $\gamma_y(\tau)$ are called  {\it Floquet multipliers}
of $(\tau,y)$. By Proposition 1.6.13 of \cite{Eke1}, the Floquet
multipliers with their multiplicities and Krein type numbers of
$(\tau,y)\in\mathcal{J}(\Sigma,\alpha)$ do not depend on the
particular choice of the Hamiltonian function in $(1.2)$. As in
Chapter 15 of \cite{Lon1}, for any symplectic matrix M, we define the elliptic height $e(M)$ of
M by the total algebraic multiplicity of all eigenvalues of M on the unit circle
${\bf U}$ in the complex plane ${\bf C}$. And for any $(\tau,y)\in\mathcal{J}(\Sigma,\alpha)$
we define $e(\tau,y)=e(\gamma_y(\tau))$, and call $(\tau,y)$ {\it elliptic} or {\it hyperbolic} if
$e(\tau,y)=2n$ or $e(\tau,y)=2$, respectively.

As in Definition 5.1.6 of \cite{Eke1}, a $C^{3}$ hypersurface $\Sigma$
bounding a compact convex set $U$, containing $0$ in its interior is
$(r,R)$-{\it pinched}, with $0<r\leq R$, if:\[|y|^{2} R^{-2}\leq
\frac{1}{2}(H_{2}^{\prime\prime}(x)y,y)\leq|y|^{2} r^{-2},
~\forall~x \in \Sigma.\]

For the existence, multiplicity and stability of closed
characteristics on convex compact hypersurfaces in ${\bf R}^{2n}$ we
refer to \cite{Rab1,Wei1,EkL1,EkL2,Gir1,EkH1,Szu1,LLZ1,LoZ1,WHL1} and references therein.
It is very interesting to consider closed characteristics on hypersurfaces with
special symmetries. \cite{Wan1, Liu2, Zha1} studied the multiplicity
of closed characteristics with symmetries on convex compact hypersurfaces without pinching conditions.
For the stability problem of closed characteristics with symmetries,
in \cite{HuS1} of 2009, Hu and Sun studied the
index theory and stability of periodic solutions in Hamiltonian systems with symmetries. As application
they studied the stability of figure-eight orbit due to Chenciner and Montgomery
in the planar three-body problems with equal masses. In \cite{Liu1}, Liu studied the stability of
symmetric closed characteristics on central symmetric compact
convex hypersurfaces under a pinching condition. In \cite{DoL2}, Dong and Long proved
that there exists at least one P-invariant closed
characteristic which possesses at least $2n-4\kappa$ 	
Floquet multipliers on the unit circle of the complex plane.
In this paper, we can obtain two such closed characteristics under a pinching condition:

{\bf Theorem 1.1.} {\it Assume $\Sigma\in\mathcal {H}_\kappa(2n)$ and $0<r \leq |x|\leq R,~\forall~x \in \Sigma$
with $\frac{R}{r}<\sqrt{\frac{5}{3}}$. Then there exist at least two geometrically distinct P-symmetric closed
characteristics which possess at least $2n-4\kappa$ 	
Floquet multipliers on the unit circle of the complex plane.}

{\bf Remark 1.2.} In the
above Theorem 1.1, let $\kappa=0$, the P-symmetric closed
characteristic is just symmetric and the P-fixed closed
characteristics vanish, so Theorem 1.1 covers Theorem 1.1 of \cite{Liu1}.

In this paper, let ${\bf N}, {\bf N}_0, {\bf Z}, {\bf Q}, {\bf R}$
and ${\bf C}$ denote the sets of natural integers, non-negative
integers, integers, rational numbers, real numbers and complex
numbers respectively. Denote by $a\cdot b$ and $|a|$ the standard
inner product and norm in ${\bf R}^{2n}$. Denote by $\langle
\cdot,\cdot \rangle$ and $\|\cdot\|$ the standard $L^2$-inner
product and $L^2$-norm. For an $S^1$-space $X$, we denote by
$X_{S^1}$ the homotopy quotient of $X$ by $S^1$, i.e.,
$X_{S^1}=S^\infty\times_{S^1}X$, where $S^\infty$ is the unit sphere
in an infinite dimensional {\it complex} Hilbert space. we define
the functions \bea[a] = \max{\{k \in {\bf Z} \mid k \leq
a\}},~~\{a\}=a-[a],~~E(a) = \min{\{k \in {\bf Z} \mid k \geq
a\}},~~\phi(a) = E(a) -[a].\nn\eea Specially, $\phi(a)=0$ if $a \in
{\bf Z}$, and $\phi(a)=1$ if $a\notin {\bf Z}$. We use $\Q$
coefficients for all homological modules.

\setcounter{equation}{0}
\section{A variational structure for P-invariant closed characteristics}
In the rest of this paper, we fix a $\Sigma\in
\mathcal
{H}_\kappa(2n)$. In this section, we review a variational structure for P-invariant closed characteristics
established in \cite{Liu2}.

As in \cite{Liu2}, we associate with U a convex function $H_{a}$. Consider
the fixed period problem\be\left\{\matrix{
\dot{x}(t)=JH_{a}^{\prime}(x(t)), \cr x(1/2)=Px(0). \cr
}\right.\ee
Then by Proposition 2.2 of \cite{Liu2}, nonzero solutions of (2.1) are in one to one
correspondence with P-symmetric closed characteristics with period
$\tau < a$ and P-fixed closed characteristics with period
$\frac{\tau}{2} < \frac{a}{2}$.
Let \bea L^2_\kappa\left(0,
\frac{1}{2}\right)&=&\{u=u_1\diamond u_2 \in L^2((0,\frac{1}{2}),
{\bf R}^{2n})\mid u_1\in L^2((0,\frac{1}{2}), {\bf
R}^{2n-2\kappa}),\nn\\&& u_2\in L^2((0,\frac{1}{2}), {\bf
R}^{2\kappa}), u(\frac{1}{2})=Pu(0),
\int_0^{\frac{1}{2}}u_2(t)dt=0\}\eea
Define a linear operator $\Pi_\kappa: L^2_\kappa\left(0,
\frac{1}{2}\right)\rightarrow L^2_\kappa\left(0, \frac{1}{2}\right)$
by\bea (\Pi_\kappa u)(t)&=&x_1(t)\diamond x_2(t),\nn\\x_1(t)&=&\int_0^t
u_1(\tau)d\tau-\frac{1}{2}\int_0^{\frac{1}{2}}u_1(\tau)d\tau,\nn\\x_2(t)&=&\int_0^t
u_2(\tau)d\tau-2\int_0^{\frac{1}{2}}dt\int_0^t u_2(\tau)d\tau,\nn\eea
for any $u=u_1\diamond u_2\in L^2_\kappa\left(0,
\frac{1}{2}\right)$.

The corresponding
Clarke-Ekeland dual action functional is defined
by\be\Psi_{a}(u)=\int_0^{\frac{1}{2}}\left(\frac{1}{2}J u \cdot
\Pi_\kappa u+G_a(-Ju)\right)dt,\ee
where $G_a$ is the
Fenchel transform of $H_{a}$ defined by $G_a(x)=\sup{\{ x\cdot
y-H_a(y) \mid {y\in {\bf R}^{2n}}\}}$. By Proposition 2.6 of \cite{Liu2},
$\Psi_{a}$ is $C^{1,1}$ on $L^2_\kappa\left(0, \frac{1}{2}\right)$ and satisfies the Palais-Smale condition.
Suppose $x$ is a
solution of $(2.1)$. Then $u =\dot{x}$ is a critical point of
$\Psi_{a}$. Conversely, suppose $u$ is a critical point of
$\Psi_{a}$. Then there exists a
unique $\xi\in{\bf R}^{2n}$ such that $\Pi_\kappa u-\xi$ is a
solution of (2.1). In particular, solutions of (2.1) are
in one to one correspondence with critical points of $\Psi_{a}$. Moreover, $\Psi_{a}(u)<0$ for every critical
point $u\neq 0$ of $\Psi_{a}$.

Suppose $u$ is a nonzero critical point
of $\Psi_a$. Then the formal Hessian of $\Psi_a$ at $u$ is defined
by\begin{eqnarray} Q_{a}(v,v)=\int_0^{\frac{1}{2}}\left(Jv\cdot
\Pi_\kappa v+G_a^{\prime\prime}(-Ju)Jv\cdot
Jv\right)dt,\end{eqnarray} which defines an orthogonal splitting
$L^2_\kappa\left(0, \frac{1}{2}\right)=E_{-}\oplus E_{0}\oplus
E_{+}$ of $L^2_\kappa\left(0, \frac{1}{2}\right)$ into negative,
zero and positive subspaces. The index of $u$ is defined by $i(u) =
dimE_{-}$ and the nullity of $u$ is defined by $\nu(u) = dimE_{0}$. cf. Definition 2.10 of \cite{Liu2}.

Note that we can identify $L^2_\kappa\left(0,
\frac{1}{2}\right)$ with the space $\{u\in L^2(\R/\Z,
{\bf R}^{2n})\mid u|_{(0,1/2)}\in L^2_\kappa(0,
\frac{1}{2}), u(t+\frac{1}{2})=Pu(t) \}$. Then we have a natural $S^{1}$-action on $L^2_\kappa\left(0,
\frac{1}{2}\right)$ defined by $\theta*u(t) = u(\theta+t)$, for all
$\theta \in S^{1}\equiv {\bf R}/{\bf Z}$ and $t \in {\bf R}$. By Lemma 2.8 of \cite{Liu2},
$\Psi_{a}$ is $S^{1}$-invariant. Hence if $u$ is a critical point of
$\Psi_{a}$, then the whole orbit $S^{1}\cdot u$ is formed by
critical points of $\Psi_{a}$. Denote by $crit(\Psi_{a})$ the set of
critical points of $\Psi_{a}$. Then $crit(\Psi_{a})$ is compact by
the Palais-Smale condition.

Recall that for a principal $U(1)$-bundle $E\rightarrow B$, the
Fadell-Rabinowitz index (cf. \cite{FaR1}) of E is defined to be
$\sup\{k \mid c_{1}(E)^{k-1}\neq 0\}$, where $c_1(E)\in H^2(B,{\bf
Q})$ is the first rational Chern class. For a $U(1)$-space, i.e., a
topological space $X$ with a $U(1)$-action, the Fadell-Rabinowitz
index is defined to be the index of the bundle $X\times
S^{\infty}\rightarrow X\times_{U(1)} S^{\infty}$, where
$S^{\infty}\rightarrow CP^{\infty}$ is the universal $U(1)$-bundle.

As on Page 199 of \cite{Eke1}, we choose some $\alpha\in (1,2)$ and associate with $U$ a convex function
$H(x)=j(x)^\alpha, \forall x\in\R^{2n}$. Consider the fixed period problem
\bea\left\{\matrix{ \dot{x}(t)=JH^{\prime}(x(t)), \cr
x(\frac{1}{2})=P x(0).~~~~\cr }\right.\nn\eea
The corresponding Clarke-Ekeland dual action functional on $L^2_\kappa\left(0,
\frac{1}{2}\right)$ is defined by
\bea\Psi(u)=\int_0^{\frac{1}{2}}\left(\frac{1}{2}J u \cdot
\Pi_\kappa u+H^*(-Ju)\right)dt, \forall u\in L^2_\kappa\left(0,
\frac{1}{2}\right),\nn\eea
where $H^*$ is the Fenchel
transform of $H$.

For any $\iota\in {\bf R}$, we denote
by\begin{eqnarray}\Psi^{\iota-}=\left\{w\in
L_\kappa^2(0,\frac{1}{2}) \mid
\Psi(w)<\iota\right\}.\end{eqnarray}
As in Section 2 of \cite{Liu2}, we define \begin{eqnarray}c_i = \inf\{\delta\in {\bf R}
\mid \hat{I}(\Psi^{\delta-})\geq i\}.\end{eqnarray}
where $\hat{I}$ is the Fadell-Rabinowitz index defined above. Then \[c_1\leq c_2\leq\cdots c_i\leq c_{i+1}\leq\cdots<0.\]
By Propositions 2.15 and 2.16 of \cite{Liu2}, we have

{\bf Proposition 2.1.} {\it Every $c_{i}$ is a critical value of
$\Psi$. If $c_{i} = c_{j}$ for some $i < j$, then there are
infinitely many geometrically distinct P-invariant closed
characteristics on $\Sigma$.}

{\bf Proposition 2.2.} {\it For every $i\in\N$, there is a critical point $u_\alpha$ of
$\Psi$ found in Proposition 2.1 such that
\begin{eqnarray}\Psi(u_\alpha)=c_{i},~~~C_{S^1,2i-2}(\Psi_{a}, S^1\cdot u)\neq0
\end{eqnarray}
where $u$ is a critical point of $\Psi_{a}$ corresponding to $u_\alpha$
in the natural sense. In particular, we
have $i(u)\leq 2(i-1) \leq i(u)+\nu(u)-1$.}

\setcounter{equation}{0}
\section{Index iteration theory for P-symmetric closed characteristics}

In this section, we review the index iteration theory for P-symmetric
closed characteristics which was studied in Section 3 of \cite{Liu2}.

Note that if $(\tau, y)$ is P-symmetric, then $((2m -1)\tau, y)$ is
P-symmetric for any $m\in {\bf N}$. Thus $((2m-1)\tau, y)$
corresponds to a critical point of $\Psi_a$ via Propositions $2.2$
and $2.6$ of \cite{Liu2}, we denote it by $u^{2m-1}$.
Recall that the action of a closed characteristic $(\tau,y)$ is
defined by (cf. P.190 of \cite{Eke1})\[
A(\tau,y)=\frac{1}{2}\int_0^{\tau}{(Jy\cdot \dot{y})}dt.\]

{\bf Lemma 3.1.}(cf. Lemma 3.1 of \cite{Liu2}) {\it Suppose $u^{2m-1}$ is a nonzero critical point
of $\Psi_{a}$ such that $u$ corresponds to P-symmetric closed
characteristic $(\tau, y)$. Let $H_{2}(x) = j^{2}(x)$, where $j$ is
the gauge function of $\Sigma$. And by $(21)$ in $P.191$ of
\cite{Eke1}, $\tau_2=A(\tau,y)$. Then $i(u^{2m-1})$ equals the index
of the following quadratic
form\begin{eqnarray}q_{{(2m-1)\tau_2/2},\kappa}(v,v):
=\int_0^{(2m-1)\tau_2/2}[(Jv,
\Pi_{{(2m-1)\tau_2/2},\kappa}v)+(H_{2}^{\prime\prime}(y(t))^{-1}Jv,
Jv)]dt .
\end{eqnarray}
where $v \in L^2_{\kappa}(0, {(2m-1)\tau_2/2})$, the definitions of
$q_{{(2m-1)\tau_2/2},\kappa}$, $\Pi_{{(2m-1)\tau_2/2},\kappa}$,
$L^2_{\kappa}(0, {(2m-1)\tau_2/2})$ are as in section 3 of
\cite{DoL2}. Moreover, we have $\nu(u^{2m-1}) = \nu
\left(q_{(2m-1)\tau_2/2,\kappa}\right)-1$.}

Now we consider the linear Hamiltonian system \be\left\{\matrix{
\dot{\xi}(t)=JA(t)\xi,~~~~~~~~ \cr A(t+\tau_2/2)=PA(t)P. \cr}
\right.\ee where $A(t)=H_2^{\prime\prime}(y(t))$. Denote by
$i_P^E(A,k)$ and $\nu_P^E(A,k)$ the index and nullity of the $k$-th
iteration of the system (3.2) defined by Dong and Long (cf.
Definition 3.4 of \cite{DoL2}). Denote by $i_{P, 1}(\gamma_A^{k,P})$
and $\nu_{P, 1}(\gamma_A^{k,P})$ the P-index and P-nullity of the
$k$-th iteration of the system $(3.2)$ defined by Dong and Long (cf.
Section 3 of \cite{DoL1}), where $\gamma_A$ is the fundamental
solution of (3.2) with $\gamma_A(0)=I_{2n}$. Then we have

{\bf Theorem 3.2.}(cf. Theorem 3.2 of \cite{Liu2}) {\it If $u^{2m-1}$ is a nonzero critical point of $\Psi_{a}$ such
that $u$ corresponds to P-symmetric closed characteristic $(\tau,
y)$. Then we have\bea &&i(u^{2m-1})=i_P^E(A,2m-1)=i_{P,
1}(\gamma_A^{2m-1,P})-\kappa, \nn\\&&\nu(u^{2m-1})=\nu_P^E(A,2m-1)-1=\nu_{P,
1}(\gamma_A^{2m-1,P})-1.\eea}

Now we compute $i(u^{2m-1})$ via the index iteration method in
\cite{Lon1} and \cite{DoL1}. First we recall briefly an index theory
for symplectic paths. All the details can be found in \cite{Lon1}, \cite{DoL1} and \cite{Liu2}.

In the following of this section, we assume $P$ is
some matrix of pattern $(-I_{2s-2t})\diamond I_{2t}$, where $0\leq
t\leq s$.

As usual, the symplectic group $Sp(2n)$ is defined by\[Sp(2n) = \{M
\in GL(2n,{\bf R}) \mid M^{T}JM = J\},\] whose topology is induced
from that of ${\bf R}^{4n^2}$. For $\tau>0$ we are interested in
paths in $Sp(2n)$:\[\mathcal {P}_{\tau}(2n) = \{\gamma \in C([0,
\tau], Sp(2n)) \mid \gamma(0) = I_{2n}\},\] which is equipped with
the topology induced from that of $Sp(2n)$. The following real
function was introduced in \cite{DoL1}:\[D_{P,\omega}(M) =
(-1)^{n-1}\bar{\omega}^n det(M-\omega P),~\forall~\omega \in {\bf
U}, M \in Sp(2n).\]where ${\bf
U}$ is the unit circle in the complex plane. Thus for any $\omega\in {\bf U}$ the following
codimension $1$ hypersurface in $Sp(2n)$ is defined in
\cite{DoL1}:\[Sp(2n)^0_{P,\omega} = \{M \in Sp(2n) \mid
D_{P,\omega}(M) = 0\}.\] For any $M \in Sp(2n)^0_{P,\omega}$, we
define a co-orientation of $Sp(2n)^0_{P,\omega}$ at $M$ by the
positive direction $\frac{d}{dt}M e^{t\epsilon J} |_{t=0}$ of the
path $M e^{t\epsilon J}$ with $0\leq t \leq 1$ and $\epsilon > 0$
being sufficiently small. Let\begin{eqnarray}&Sp(2n)_{P,\omega}^*¥ø
= Sp(2n)\setminus Sp(2n)^0_{P,\omega} ,&\nn\\&\mathcal
{P}^*_{P,\tau,\omega}(2n) = \{\gamma \in \mathcal {P}_{\tau} (2n)
\mid \gamma_{\tau}\in Sp(2n)^*_{P,\omega}¥ø\},&\nn\\&\mathcal
{P}^0_{P,\tau,\omega}(2n)=\mathcal {P}_{\tau} (2n)\setminus \mathcal
{P}^*_{P,\tau,\omega}(2n).&\nn\end{eqnarray} For any two continuous
arcs $\xi$ and $\eta : [0, \tau ] \longrightarrow Sp(2n)$ with
$\xi(\tau) = \eta(0)$, it is defined as
usual:\\$${\eta}*{\xi}(t)=\left\{ \begin{array}{ll}\xi(2t),~if~0\leq t\leq \tau/2,\\
\eta(2t-1),~if~\tau/2\leq t \leq \tau.\end{array}\right.$$ Given any
two $2m_k \times 2m_k$ matrices of square block form $M_{k}=\left(
              \begin{array}{cc}
                A_k & B_k \\
                C_k & D_k\\
              \end{array}
            \right)$
with $k = 1, 2$, as in \cite{Lon1}, the $\diamond$-product of $M_1$
and $M_2$ is defined by the following $2(m_1 +m_2) \times 2(m_1
+m_2)$ matrix $M_1 \diamond M_2$:\\$$M_1 \diamond M_2=\left(
                            \begin{array}{cccc}
                              A_1& 0 & B_1 & 0 \\
                              0 & A_2 & 0 & B_2 \\
                              C_1 & 0 & D_1 & 0 \\
                              0 & C_2 & 0 & D_2 \\
                            \end{array}
                          \right)$$
Denote by $M^{\diamond k}$ the $k$-fold $\diamond$-product $M
\diamond \cdots \diamond M$. Note that the $\diamond$-product of any
two symplectic matrices is symplectic. For any two paths $\gamma_{j}
\in \mathcal {P}_{\tau}(2n_j)$ with $j = 0$ and $1$, let $\gamma_1
\diamond \gamma_2(t)=\gamma_1(t)\diamond \gamma_2(t)$ for all $t \in
[0, \tau]$.

A special path $\xi_n \in \mathcal {P}_{\tau} (2n)$ is defined
by\\$$\xi_n(t)=\left(
                 \begin{array}{cc}
                   2-\frac{t}{\tau}& 0 \\
                   0 & (2-\frac{t}{\tau})^{-1} \\
                 \end{array}
               \right)^{\diamond n},~for~~0 \leq t \leq \tau.$$
{\bf Definition 3.3.}(cf. \cite{DoL1}, also Definition 3.3 of \cite{Liu2})
For any $\omega \in {\bf U}$
and $M \in Sp(2n)$, via Definition 5.4.4 in \cite{Lon1}, we
define\begin{eqnarray}\nu_{P,\omega}(M) = dim_{{\bf C}} ker_{{\bf
C}}(M - \omega P).\end{eqnarray} For any $\tau
> 0$ and $\gamma \in \mathcal {P}_{\tau} (2n)$,
define\begin{eqnarray}\nu_{P,\omega}(\gamma) =\nu_{\omega}(\gamma
P)= \nu_{P,\omega}(\gamma(\tau )).\end{eqnarray} If $\gamma \in
\mathcal {P}^*_{P,\tau,\omega}(2n)$,
define\begin{eqnarray}i_{P,\omega}(\gamma) = [Sp(2n)^0_{P,\omega}:
\gamma\ast \xi_n],\end{eqnarray} where the right hand side of
$(3.6)$ is the usual homotopy intersection number, and the
orientation of $\gamma\ast \xi_n$ is its positive time direction
under homotopy with fixed end points.

If $\gamma\in \mathcal {P}^0_{P,\tau,\omega}(2n)$, we let $\mathcal
{F}(\gamma)$ be the set of all open neighborhoods of $\gamma$ in
$\mathcal {P}_{\tau}(2n)$, and define
\begin{eqnarray}
i_{P,\omega}(\gamma) = \sup_{ U\in \mathcal {F}(\gamma)}{
\inf{\{i_{P,\omega}(\beta) \mid \beta \in U\cap\mathcal
{P}^*_{P,\tau,\omega}(2n)\}}}.
\end{eqnarray}
Then\begin{eqnarray}(i_{P,\omega}(\gamma), \nu_{P,\omega}(\gamma))
\in {\bf Z} \times \{0, 1, \cdots , 2n\},\end{eqnarray} is called
the P-index function of $\gamma$ at $\omega$.

For any $M \in Sp(2n)$ and $\omega \in {\bf U}$, {\it the splitting
numbers} $S^{\pm}_M (P,\omega)$ of $M$ at $(P,\omega)$ are defined
by
\begin{eqnarray}
S^{\pm}_M (P,\omega)=\lim_{\epsilon\rightarrow 0^+}{i_{P,\omega
\exp{(\pm\sqrt{-1}\epsilon)}}(\gamma)-i_{P,\omega}(\gamma)},
\end{eqnarray}
for any path $\gamma \in \mathcal {P}_{\tau} (2n)$ satisfying
$\gamma(\tau) = M$.

Let $\Omega^0(M)$ be the path connected component containing $M =
\gamma(\tau)$ of the set\begin{eqnarray} \Omega(M) = \{N \in Sp(2n)
\mid  \sigma(N)
\cap {\bf U} = \sigma(M) \cap {\bf U}~and~~~~\\
\nu_{\lambda}(N)=\nu_{\lambda}(M),\forall \lambda \in \sigma(M) \cap
{\bf U}\}\end{eqnarray} Here $\Omega^0(M)$ is called the {\it
homotopy component} of $M$ in $Sp(2n)$.

In \cite{Lon1}, the following symplectic matrices were introduced as
basic normal forms: \begin{eqnarray}D(\lambda)=\left(
                                                    \begin{array}{cc}
                                                      \lambda & 0 \\
                                                      0 & \lambda^{-1} \\
                                                    \end{array}
                                                  \right),
                                                  \lambda=\pm 2,\\
N_1(\lambda,b)=\left(
 \begin{array}{cc}
  \lambda & b \\
   0& \lambda \\
   \end{array}
   \right),\lambda=\pm1,b=\pm1,0,\\
R(\theta)=\left(
            \begin{array}{cc}
              \cos{\theta} & -\sin{\theta} \\
              \sin{\theta} & \cos{\theta} \\
            \end{array}
          \right),\theta\in (0,\pi)\cup(\pi,2\pi),\\
N_2(\omega,B)=\left(
                \begin{array}{cc}
                  R(\theta) & B \\
                  0 &  R(\theta) \\
                \end{array}
              \right),\theta\in (0,\pi)\cup(\pi,2\pi),
\end{eqnarray}
where $B=\left(
           \begin{array}{cc}
             b_1 & b_2 \\
             b_3 & b_4 \\
           \end{array}
         \right)
$ with $b_i\in {\bf R}$ and $b_2\neq b_3$.

Splitting numbers possess the following properties:

{\bf Lemma 3.4.}(cf. Proposition 3.8 of \cite{DoL1}) {\it Let
$(p_\omega(MP), q_\omega(MP))$ denote the Krein type of $MP$ at
$\omega$. For any $M\in Sp(2n)$ and $\omega\in {\bf U}$, the
splitting numbers $S^{\pm}_M (P,\omega)$ are well defined and
satisfy the following properties.

(i) $S^{\pm}_M (P,\omega)=S^{\pm}_{MP}(\omega)$, where the
right-hand side is the splitting numbers given by Definition 9.1.4
of \cite{Lon1}.

(ii) $S^{+}_M (P,\omega)-S^{-}_M (P,\omega)=p_\omega(MP)-
q_\omega(MP)$.

(iii) $S^{\pm}_M (P,\omega)=S^{\pm}_N (P,\omega)$ if $NP\in
\Omega^0(MP)$.

(iv) $S^{\pm}_{M_1\diamond
M_2}(P,\omega)=S^{\pm}_{M_1}(P_1,\omega)+S^{\pm}_{M_2}(P_2,\omega)$
for $M_j, P_j\in Sp(2n_j)$ with $n_j\in\{1,\cdots,n\}$ satisfying
$P=P_1\diamond P_2$ and $n_1+n_2=n$.

(v) $S^{\pm}_M (P,\omega)=0$ if $\omega\notin \sigma(MP)$.}

We have the following

{\bf Lemma 3.5.}(cf. Theorem 1.8.10 of \cite{Lon1}) {\it For any $M
\in Sp(2n)$, there is a path $f : [0, 1] \rightarrow \Omega^0(M)$
such that $f(0) = M$ and\begin{eqnarray} f(1)=M_1\diamond \cdots
\diamond M_l,
\end{eqnarray}
where each $M_i$ is a basic normal form listed in (3.12)-(3.15) for
$1 \leq i \leq l$. \hfill\hb}

By Proposition 3.10 of \cite{DoL1}, we have the Bott-type formula for $(P,\omega)$-index:

{\bf Lemma 3.6.} {\it For any $\gamma \in \mathcal {P}_{\tau} (2n)$, $z \in {\bf U}$
and $m\in\N$, we have\bea i_{P^m,z}(\gamma^{m,P})&=&\sum_{\omega^{m}=z}
i_{P,\omega}(\gamma),\nn\\ \nu_{P^m,z}(\gamma^{m,P})&=&\sum_{\omega^{m}=z}
\nu_{P,\omega}(\gamma).\nn\eea}

Now we deduce the index iteration formula for each case in
(3.12)-(3.15). Note that the splitting numbers are computed in List
9.1.12 of \cite{Lon1}. Let $M=\gamma(\tau)$.

{\bf Case 1.} $MP$ is conjugate to a matrix $\left(
              \begin{array}{cc}
                1 & b \\
                0 & 1\\
              \end{array}
            \right)$ for some $b>0$.

In this case, we have $(S^{+}_{M}(P,1),S^{-}_{M}(P,1))=(1,1)$ by List
9.1.12 of \cite{Lon1} and Lemma 3.4 (i), (iii). Thus
by Lemma 3.4 (v) and Lemma 3.6, we
have\begin{eqnarray}&i_{P,1}(\gamma^{2m-1,P})&=\sum_{\omega^{2m-1}=1}
i_{P,\omega}(\gamma)=\sum_{k=0}^{2m-2}i_{P,e^{2k\pi i/(2m-1)}}(\gamma)=(2m-1)(i_{P,1}(\gamma)+1)-1,
\nn\\&\nu_{P,1}(\gamma^{2m-1,P})&=1.\end{eqnarray}

{\bf Case 2.} $MP = I_{2}$, the $2\times2$ identity matrix.

In this case, we have $(S^{+}_{M}(P,1),S^{-}_{M}(P,1))=(1,1)$ by List
9.1.12 of \cite{Lon1} and Lemma 3.4 (i), (iii). Thus
as in Case 1, we
have\begin{eqnarray}i_{P,1}(\gamma^{2m-1,P})=(2m-1)(i_{P,1}(\gamma)+1)-1,~~~~~
\nu_{P,1}(\gamma^{2m-1,P})=2.
\end{eqnarray}

{\bf Case 3.} $M P$ is conjugate to a matrix  $\left(
              \begin{array}{cc}
                1 & b \\
                0 & 1\\
              \end{array}
            \right)$ for some $b<0$.

In this case, we have $(S^{+}_{M}(P,1),S^{-}_{M}(P,1))=(0,0)$ by List
9.1.12 of \cite{Lon1} and Lemma 3.4 (i), (iii). Thus
by Lemma 3.4 (v) and Lemma 3.6, we have
\begin{eqnarray}&i_{P,1}(\gamma^{2m-1,P})&=\sum_{\omega^{2m-1}=1}
i_{P,\omega}(\gamma)=\sum_{k=0}^{2m-2}i_{P,e^{2k\pi i/(2m-1)}}(\gamma)=(2m-1)i_{P,1}(\gamma),\nn
\\&\nu_{P,1}(\gamma^{2m-1,P})&=1.\end{eqnarray}

{\bf Case 4.} $MP$ is conjugate to a matrix  $\left(
              \begin{array}{cc}
                -1 & b \\
                0 & -1\\
              \end{array}
            \right)$ for some $b<0$.

In this case, we have $(S^{+}_{M}(P,-1),S^{-}_{M}(P,-1))=(1,1)$ by List
9.1.12 of \cite{Lon1} and Lemma 3.4 (i), (iii).
Thus by Lemma 3.4 (v) and Lemma 3.6, we
have\begin{eqnarray}&i_{P,1}(\gamma^{2m-1,P})&=\sum_{\omega^{2m-1}=1}
i_{P,\omega}(\gamma)=\sum_{k=0}^{2m-2}i_{P,e^{2k\pi i/(2m-1)}}(\gamma)=(2m-1)i_{P,1}(\gamma),\nn\\&
\nu_{P,1}(\gamma^{2m-1,P})&=0.\end{eqnarray}

{\bf Case 5.} $MP = -I_{2}$.

In this case, we have $(S^{+}_{M}(P,-1),S^{-}_{M}(P,-1))=(1,1)$ by List
9.1.12 of \cite{Lon1} and Lemma 3.4 (i), (iii).
Thus as in Case 4, we
have\begin{eqnarray}i_{P,1}(\gamma^{2m-1,P})=(2m-1)i_{P,1}(\gamma),~~~~~~~\nu_{P,1}(\gamma^{2m-1,P})=0.\end{eqnarray}

{\bf Case 6.} $MP$ is conjugate to a matrix  $\left(
              \begin{array}{cc}
                -1 & b \\
                0 & -1\\
              \end{array}
            \right)$ for some $b>0$.

In this case, we have $(S^{+}_{M}(P,-1),S^{-}_{M}(P,-1))=(0,0)$ by List
9.1.12 of \cite{Lon1} and Lemma 3.4 (i), (iii).
Thus by Lemma 3.4 (v) and Lemma 3.6, we
have\begin{eqnarray}&i_{P,1}(\gamma^{2m-1,P})&=\sum_{\omega^{2m-1}=1}
i_{P,\omega}(\gamma)=\sum_{k=0}^{2m-2}i_{P,e^{2k\pi i/(2m-1)}}(\gamma)=(2m-1)i_{P,1}(\gamma),
\\&\nu_{P,1}(\gamma^{2m-1,P})&=0.\end{eqnarray}

{\bf Case 7.} $MP =\left(
              \begin{array}{cc}
                \cos\theta & -\sin\theta\\
                 \sin\theta& \cos\theta\\
              \end{array}
            \right)$ with some $\theta \in (0, \pi) \cup (\pi, 2\pi)$.

In this case, we have
$(S^{+}_{M}(P,e^{\sqrt{-1}\theta}),S^{-}_{M}(P,e^{\sqrt{-1}\theta}))=(0,1)$ by List
9.1.12 of \cite{Lon1} and Lemma 3.4 (i), (iii).
Thus by Lemma 3.4 (v) and Lemma 3.6, we
have\begin{eqnarray}i_{P,1}(\gamma^{2m-1,P})&=&\sum_{\omega^{2m-1}=1}
i_{P,\omega}(\gamma)=\sum_{k=1}^{2m-1}i_{P,e^{2k\pi i/(2m-1)}}(\gamma)\nn\\&=&
\sum_{0<2k<\frac{(2m-1)\theta}{\pi}}i_{P,1}(\gamma)+\sum_{\frac{(2m-1)\theta}{\pi}\leq
2k\leq
\frac{(2m-1)(2\pi-\theta)}{\pi}}(i_{P,1}(\gamma)-1)\nn\\&+&\sum_{\frac{(2m-1)(2\pi-\theta)}{\pi}<
2k\leq
4m-2}i_{P,1}(\gamma)\nn\\&=&(2m-1)(i_{P,1}(\gamma)-1)+2E(\frac{(2m-1)\theta}{2\pi})-1,\nn\\
\nu_{P,1}(\gamma^{2m-1,P})&=&2-2\phi(\frac{(2m-1)\theta}{2\pi}),
\end{eqnarray}
where the function $E(\cdot)$ is defined as in Section 1.

Provided $\theta\in(0,\pi)$. When $\theta\in (\pi,2\pi)$, we have
\begin{eqnarray}i_{P,1}(\gamma^{2m-1,P})&=&\sum_{\omega^{2m-1}=1}
i_{P,\omega}(\gamma)=\sum_{k=1}^{2m-1}i_{P,e^{2k\pi i/(2m-1)}}(\gamma)\nn\\&=&
\sum_{0<2k\leq\frac{(2m-1)(2\pi-\theta)}{\pi}}i_{P,1}(\gamma)+\sum_{\frac{(2m-1)(2\pi-\theta)}{\pi}<
2k <
\frac{(2m-1)\theta}{\pi}}(i_{P,1}(\gamma)+1)\nn\\&+&\sum_{\frac{(2m-1)\theta}{\pi}\leq
2k\leq
4m-2}i_{P,1}(\gamma)\nn\\&=&(2m-1)(i_{P,1}(\gamma)-1)+2E(\frac{(2m-1)\theta}{2\pi})-1,\nn\\
\nu_{P,1}(\gamma^{2m-1,P})&=&2-2\phi(\frac{(2m-1)\theta}{2\pi}).
\end{eqnarray}

{\bf Case 8.} $MP =\left(
              \begin{array}{cc}
                R(\theta) & B\\
                 0 & R(\theta)\\
              \end{array}
            \right)$ with some $\theta \in (0, \pi) \cup (\pi,
            2\pi)$ and $B=\left(
              \begin{array}{cc}
                b_{1} & b_{2}\\
                 b_{3} & b_{4}\\
              \end{array}
            \right) \in {\bf R}^{2\times 2}$, such
that $(b_{2}- b_{3}) \sin\theta < 0$.

In this case, we have
$(S^{+}_{M}(P,e^{\sqrt{-1}\theta}),S^{-}_{M}(P,e^{\sqrt{-1}\theta}))=(1,1)$ by List
9.1.12 of \cite{Lon1} and Lemma 3.4 (i), (iii).
Thus by Lemma 3.4 (v) and Lemma 3.6, we
have\begin{eqnarray}i_{P,1}(\gamma^{2m-1,P})&=&\sum_{\omega^{2m-1}=1}
i_{P,\omega}(\gamma)=\sum_{k=1}^{2m-1}i_{P,e^{2k\pi i/(2m-1)}}(\gamma)\nn\\
&=&(2m-1)i_{P,1}(\gamma)+2\phi(\frac{(2m-1)\theta}{2\pi})-2,\nn\\
\nu_{P,1}(\gamma^{2m-1,P})&=&2-2\phi(\frac{(2m-1)\theta}{2\pi}).\end{eqnarray}
Here from the first line to the second line of Equation (3.26), we used that
if $\frac{(2m-1)\theta}{2\pi}\notin\N$, then
$i_{P,e^{2k\pi i/(2m-1)}}(\gamma)=i_{P,1}(\gamma)$ for all $1\leq k \leq 2m-1$ and
if $\frac{(2m-1)\theta}{2\pi}\in\N$,
then exactly two of the $i_{P,e^{2k\pi i/(2m-1)}}(\gamma)$'s equal to $i_{P,1}(\gamma)-1$
and the other ones equal to $i_{P,1}(\gamma)$.

{\bf Case 9.} $MP =\left(
              \begin{array}{cc}
                R(\theta) & B\\
                 0 & R(\theta)\\
              \end{array}
            \right)$ with some $\theta \in (0, \pi) \cup (\pi,
            2\pi)$ and $B=\left(
              \begin{array}{cc}
                b_{1} & b_{2}\\
                 b_{3} & b_{4}\\
              \end{array}
            \right) \in \textbf{R}^{2\times 2}$, such
that $(b_{2}- b_{3}) \sin\theta > 0$.

In this case, we have
$(S^{+}_{M}(P,e^{\sqrt{-1}\theta}),S^{-}_{M}(P,e^{\sqrt{-1}\theta}))=(0,0)$ by List
9.1.12 of \cite{Lon1} and Lemma 3.4 (i), (iii).
Thus by Lemma 3.4 (v) and Lemma 3.6, we
have\begin{eqnarray}&i_{P,1}(\gamma^{2m-1,P})&=\sum_{\omega^{2m-1}=1}
i_{P,\omega}(\gamma)=\sum_{k=1}^{2m-1}i_{P,e^{2k\pi i/(2m-1)}}(\gamma)
=(2m-1)i_{P,1}(\gamma),\nn\\
&\nu_{P,1}(\gamma^{2m-1,P})&=2-2\phi(\frac{(2m-1)\theta}{2\pi}),\end{eqnarray}

{\bf Case 10.} $MP$ is hyperbolic, i.e., $\textbf{U}\cap
\sigma(MP)=\emptyset$.

In this case, by Lemma 3.4 (v) and Lemma 3.6, we
have\begin{eqnarray}&i_{P,1}(\gamma^{2m-1,P})&=\sum_{\omega^{2m-1}=1}
i_{P,\omega}(\gamma)=\sum_{k=1}^{2m-1}i_{P,e^{2k\pi i/(2m-1)}}(\gamma)
=(2m-1)i_{P,1}(\gamma),\nn\\
&\nu_{P,1}(\gamma^{2m-1,P})&=0.\end{eqnarray}

The following new theorem will play a crucial role in our proof of
Theorem 1.1:

{\bf Theorem 3.7.} {\it Suppose $u^{2m-1}$ is a nonzero critical
point of $\Psi_{a}$ such that $u$ corresponds to a P-symmetric closed characteristic $(\tau, y)$. Then
we have $i(u^{3})-3i(u)\leq 2\kappa+2n$ and $i(u^{3})+\nu(u^{3})-3(i(u)+\nu(u))\geq 2\kappa+2-2n$.
In particular, we have the following

(i) if $i(u^{3})-3i(u)\geq 2n$, then $e(\tau, y)\geq 2n-2\kappa$.

(ii) if $i(u^{3})+\nu(u^{3})-3(i(u)+\nu(u))\leq 6\kappa+2-2n$, then $e(\tau, y)\geq 2n-4\kappa$.}

{\bf Proof.} Firstly, we compute $i_{P,1}(\gamma^{3,P})-3i_{P,1}(\gamma)$ and
$i_{P,1}(\gamma^{3,P})+\nu_{P,1}(\gamma^{3,P})-3(i_{P,1}(\gamma)+\nu_{P,1}(\gamma))$
for any symplectic path $\gamma\in \mathcal {P}_{\tau} (2n)$
satisfying $\gamma(\tau)=M$. We consider each of the above cases.

{\bf Case 1.}  $MP$ is conjugate to a matrix $\left(
              \begin{array}{cc}
                1 & b \\
                0 & 1\\
              \end{array}
            \right)$ for some $b>0$.

In this case, we have \bea
&&i_{P,1}(\gamma^{3,P})-3i_{P,1}(\gamma)=3(i_{P,1}(\gamma)+1)-1-3i_{P,1}(\gamma)=2,\nn\\
&&i_{P,1}(\gamma^{3,P})+\nu_{P,1}(\gamma^{3,P})-3(i_{P,1}(\gamma)+\nu_{P,1}(\gamma))=
3(i_{P,1}(\gamma)+1)-3(i_{P,1}(\gamma)+1)=0.\nn\eea

{\bf Case 2.} $MP = I_{2}$, the $2\times2$ identity matrix.

In this case, we have  \bea
&&i_{P,1}(\gamma^{3,P})-3i_{P,1}(\gamma)=3(i_{P,1}(\gamma)+1)-1-3i_{P,1}(\gamma)=2,\nn\\
&&i_{P,1}(\gamma^{3,P})+\nu_{P,1}(\gamma^{3,P})-3(i_{P,1}(\gamma)+\nu_{P,1}(\gamma))=3(i_{P,1}(\gamma)+1)+1-3(i_{P,1}(\gamma)+2)
=-2.\nn\eea

{\bf Case 3.} $M P$ is conjugate to a matrix  $\left(
              \begin{array}{cc}
                1 & b \\
                0 & 1\\
              \end{array}
            \right)$ for some $b<0$.

In this case, we have  \bea
&&i_{P,1}(\gamma^{3,P})-3i_{P,1}(\gamma)=3i_{P,1}(\gamma)-3i_{P,1}(\gamma)=0,\nn\\
&&i_{P,1}(\gamma^{3,P})+\nu_{P,1}(\gamma^{3,P})-3(i_{P,1}(\gamma)+\nu_{P,1}(\gamma))=3i_{P,1}(\gamma)+1-3(i_{P,1}(\gamma)+1)
=-2.\nn\eea

{\bf Case 4.} $MP$ is conjugate to a matrix  $\left(
              \begin{array}{cc}
                -1 & b \\
                0 & -1\\
              \end{array}
            \right)$ for some $b<0$.

In this case, we have \bea
&&i_{P,1}(\gamma^{3,P})-3i_{P,1}(\gamma)=3i_{P,1}(\gamma)-3i_{P,1}(\gamma)=0,\nn\\
&&i_{P,1}(\gamma^{3,P})+\nu_{P,1}(\gamma^{3,P})-3(i_{P,1}(\gamma)+\nu_{P,1}(\gamma))=3i_{P,1}(\gamma)-3i_{P,1}(\gamma)
=0.\nn\eea

{\bf Case 5.} $MP = -I_{2}$.

In this case, we have \bea
&&i_{P,1}(\gamma^{3,P})-3i_{P,1}(\gamma)=3i_{P,1}(\gamma)-3i_{P,1}(\gamma)=0,\nn\\
&&i_{P,1}(\gamma^{3,P})+\nu_{P,1}(\gamma^{3,P})-3(i_{P,1}(\gamma)+\nu_{P,1}(\gamma))=3i_{P,1}(\gamma)-3i_{P,1}(\gamma)
=0.\nn\eea

{\bf Case 6.} $MP$ is conjugate to a matrix  $\left(
              \begin{array}{cc}
                -1 & b \\
                0 & -1\\
              \end{array}
            \right)$ for some $b>0$.

In this case, we have \bea
&&i_{P,1}(\gamma^{3,P})-3i_{P,1}(\gamma)=3i_{P,1}(\gamma)-3i_{P,1}(\gamma)=0,\nn\\
&&i_{P,1}(\gamma^{3,P})+\nu_{P,1}(\gamma^{3,P})-3(i_{P,1}(\gamma)+\nu_{P,1}(\gamma))=3i_{P,1}(\gamma)-3i_{P,1}(\gamma)
=0.\nn\eea

{\bf Case 7.} $MP =\left(
              \begin{array}{cc}
                \cos\theta & -\sin\theta\\
                 \sin\theta& \cos\theta\\
              \end{array}
            \right)$ with some $\theta \in (0, \pi) \cup (\pi, 2\pi)$.

In this case, we have  \bea
&&i_{P,1}(\gamma^{3,P})-3i_{P,1}(\gamma)=3(i_{P,1}(\gamma)-1)+2E(\frac{3\theta}{2\pi})-1-3i_{P,1}(\gamma)\leq2,\nn\\
&&i_{P,1}(\gamma^{3,P})+\nu_{P,1}(\gamma^{3,P})-3(i_{P,1}(\gamma)+\nu_{P,1}(\gamma))\nn\\&&=\left(3(i_{P,1}(\gamma)-1)+
2E(\frac{3\theta}{2\pi})-1+2-2\phi(\frac{3\theta}{2\pi})\right)-3i_{P,1}(\gamma)\nn\\
&&=-2+2E(\frac{3\theta}{2\pi})-2\phi(\frac{3\theta}{2\pi})\geq-2.\nn\eea

{\bf Case 8.} $MP =\left(
              \begin{array}{cc}
                R(\theta) & B\\
                 0 & R(\theta)\\
              \end{array}
            \right)$ with some $\theta \in (0, \pi) \cup (\pi,
            2\pi)$ and $B=\left(
              \begin{array}{cc}
                b_{1} & b_{2}\\
                 b_{3} & b_{4}\\
              \end{array}
            \right) \in {\bf R}^{2\times 2}$, such
that $(b_{2}- b_{3}) \sin\theta < 0$.

In this case, we have\bea
&&i_{P,1}(\gamma^{3,P})-3i_{P,1}(\gamma)=3i_{P,1}(\gamma)+2\phi(\frac{3\theta}{2\pi})-2-3i_{P,1}(\gamma)\leq0,\nn\\
&&i_{P,1}(\gamma^{3,P})+\nu_{P,1}(\gamma^{3,P})-3(i_{P,1}(\gamma)+\nu_{P,1}(\gamma))=3i_{P,1}(\gamma)-3i_{P,1}(\gamma)
=0.\nn\eea

{\bf Case 9.} $MP =\left(
              \begin{array}{cc}
                R(\theta) & B\\
                 0 & R(\theta)\\
              \end{array}
            \right)$ with some $\theta \in (0, \pi) \cup (\pi,
            2\pi)$ and $B=\left(
              \begin{array}{cc}
                b_{1} & b_{2}\\
                 b_{3} & b_{4}\\
              \end{array}
            \right) \in {\bf R}^{2\times 2}$, such
that $(b_{2}- b_{3}) \sin\theta > 0$.

In this case, we have\bea
&&i_{P,1}(\gamma^{3,P})-3i_{P,1}(\gamma)=3i_{P,1}(\gamma)-3i_{P,1}(\gamma)=0,\nn\\
&&i_{P,1}(\gamma^{3,P})+\nu_{P,1}(\gamma^{3,P})-3(i_{P,1}(\gamma)+\nu_{P,1}(\gamma))=3i_{P,1}(\gamma)+2-2\phi(\frac{3\theta}{2\pi})-3i_{P,1}(\gamma)
\geq 0.\nn\eea

{\bf Case 10.} $MP$ is hyperbolic, i.e., ${\bf U}\cap
\sigma(MP)=\emptyset$.

In this case, we have\bea
&&i_{P,1}(\gamma^{3,P})-3i_{P,1}(\gamma)=3i_{P,1}(\gamma)-3i_{P,1}(\gamma)=0,\nn\\
&&i_{P,1}(\gamma^{3,P})+\nu_{P,1}(\gamma^{3,P})-3(i_{P,1}(\gamma)+\nu_{P,1}(\gamma))=3i_{P,1}(\gamma)-3i_{P,1}(\gamma)
=0.\nn\eea

Now for a P-symmetric closed characteristic $(\tau, y)$, $H_2, \tau_2$ is defined as in
Lemma 3.1. Its fundamental solution $\gamma=\gamma_y:
[0,\tau_2/2]\rightarrow Sp(2n)$ with $\gamma_y(0)=I_{2n}$ of the
linearized Hamiltonian
system\begin{eqnarray}\dot{\gamma_y}(t)=JH_2^{\prime\prime}(y(t))\gamma_y(t),~\forall~t\in
{\bf R},\end{eqnarray} is called the {\it associated symplectic
path} of P-symmetric closed characteristic $(\tau, y)$.
By Lemma 3.5, we suppose $\gamma(\frac{\tau_2}{2})P$ is conjugate to
$N_{1}(1,1)^{\diamond p_{-}}\diamond N_{1}(1,-1)^{\diamond
p_{+}}\diamond(I_{2p_{0}})\diamond R(\theta_{1})\diamond\cdots
\diamond R(\theta_{r})\diamond
N_2(\omega_1,B_1)\diamond\cdots\diamond N_2(\omega_s,B_s)\diamond
M_0$ with $\sigma(M_0)\cap ({\bf U}-\{-1\})=\emptyset$. Combining the
above cases 1-10, we obtain\bea &&i_{P,1}(\gamma^{3,P})-3i_{P,1}(\gamma)\leq 2p_{-}+2p_{0}+2r\leq 2n,\nn\\
&&i_{P,1}(\gamma^{3,P})+\nu_{P,1}(\gamma^{3,P})-3(i_{P,1}(\gamma)+\nu_{P,1}(\gamma))\geq
-2p_{0}-2p_{+}-2r\geq -2n.\eea
And by Theorem 3.2, we have
$\nu(u^3)=\nu_{P,1}(\gamma^{3,P})-1$, $\nu(u)=\nu_{P,1}(\gamma)-1$,
$i(u^3)=i_{P,1}(\gamma^{3,P})-\kappa$, $i(u)=i_{P,1}(\gamma)-\kappa$. Then we
obtain\bea i(u^{3})-3i(u)&=&2\kappa+i_{P,1}(\gamma^{3,P})-3i_{P,1}(\gamma)\nn\\
&\leq&2\kappa+2p_{-}+2p_{0}+2r\leq 2\kappa+2n,\eea
and\begin{eqnarray}&&i(u^{3})+\nu(u^{3})-3(i(u)+\nu(u))\nn\\&&
=2\kappa+2+i_{P,1}(\gamma^{3,P})+\nu_{P,1}(\gamma^{3,P})-3(i_{P,1}(\gamma)+\nu_{P,1}(\gamma))\nn\\
&&\geq 2\kappa+2-2p_{0}-2p_{+}-2r
\geq 2\kappa+2-2n.\end{eqnarray}
If $i(u^{3})-3i(u)\geq 2n$, then by (3.31), we have \bea 2p_{-}+2p_{0}+2r\geq 2n-2\kappa,\eea
i.e. $e(\gamma(\frac{\tau_2}{2})P)\geq 2n-2\kappa$.

If $i(u^{3})+\nu(u^{3})-3(i(u)+\nu(u))\leq 6\kappa+2-2n$, then by (3.32), we have \bea 2p_{0}+2p_{+}+2r\geq 2n-4\kappa,\eea
i.e. $e(\gamma(\frac{\tau_2}{2})P)\geq 2n-4\kappa$. Noticing that $e(\tau, y)=e(\gamma(\tau_2))=e((\gamma(\frac{\tau_2}{2})P)^2)$,
we complete our proof.

\setcounter{equation}{0}
\section{ Proof of the main theorem}

In this section, we give the proof of the main theorem.

Firstly, we point out a minor error in Example 6 on Page 278 of \cite{DoL2} which is
useful for us:

{\bf Lemma 4.1.} {\it For any $c>0$, we have\bea
i_P^E(cI_{2n}|_{[0, s]})=2\kappa \left(E(\frac{cs}{2\pi})-1\right)+2(n-\kappa)\left(E(\frac{cs+\pi}{2\pi})-1\right).\eea}

{\bf Proof.} Note that the definition of $E(a)$ in our paper is
different from that in \cite{DoL2} by $1$. In the proof of Example 6 of \cite{DoL2}, $t=s_k=(2k\pi+\frac{3}{2}\pi)/c$
should be changed into $t=s_k=(2k\pi+\pi)/c$, and the expression
in (3.14) of \cite{DoL2} is wrong. One can easily verify our expression in (4.1) is right.

Note that since $H_2(\cdot)$ is positive homogeneous of degree-two, by the
$(r,R)$-pinched condition we have\begin{eqnarray}|x|^2R^{-2}\leq
H_{2}(x)\leq |x|^2r^{-2},~\forall~x \in
\Sigma\end{eqnarray}Comparing with the theorem of Morse-Schoenberg
in the study of geodesics, we have the following

{\bf Proposition 4.2.} {\it Let $\Sigma\in\mathcal {H}_\kappa(2n)$ which is $(r,R)$-pinched.
Suppose $u^{2m-1}$ is a nonzero critical point of $\Psi_{a}$ such
that $u$ corresponds to a P-symmetric closed characteristic $(\tau, y)$,
and $\tau_2=A(\tau,y)$ as in Lemma $3.1$. Then we have the
following\begin{eqnarray}i(u^{2m-1})\geq 2nl,~~if
~\frac{(2m-1)\tau_2}{2}>l \pi
R^{2};~~~~~~~~~~~~~~~\\
i(u^{2m-1})+\nu(u^{2m-1})\leq 2n(l-1)-1,~~if
~\frac{(2m-1)\tau_2}{2}<(l-\frac{1}{2})\pi r^{2},\end{eqnarray}
for some $l\in \N$.}

{\bf Proof.} Consider the following three quadratic forms
on $L^2_{\kappa}(0, {(2m-1)\tau_2/2})$
\bea
&&q^R_{{(2m-1)\tau_2/2},\kappa}(v,v):
=\int_0^{(2m-1)\tau_2/2}[(Jv,
\Pi_{{(2m-1)\tau_2/2},\kappa}v)+(\frac{R^2}{2}Jv,
Jv)]dt\nn\\&&q_{{(2m-1)\tau_2/2},\kappa}(v,v):
=\int_0^{(2m-1)\tau_2/2}[(Jv,
\Pi_{{(2m-1)\tau_2/2},\kappa}v)+(H_{2}^{\prime\prime}(y(t))^{-1}Jv,
Jv)]dt\nn\\&&
q^r_{{(2m-1)\tau_2/2},\kappa}(v,v):
=\int_0^{(2m-1)\tau_2/2}[(Jv,
\Pi_{{(2m-1)\tau_2/2},\kappa}v)+(\frac{r^2}{2}Jv,
Jv)]dt \nn
\eea By the $(r,R)$-pinched condition, we have
\bea q^R_{{(2m-1)\tau_2/2},\kappa}(v,v) \geq q_{{(2m-1)\tau_2/2},\kappa}(v,v)
\geq q^r_{{(2m-1)\tau_2/2},\kappa}(v,v).\nn\eea Thus we have
$i^{R}_{(2m-1)\tau_2/2}\leq i_{(2m-1)\tau_2/2} \leq
i^{r}_{(2m-1)\tau_2/2}$, where $i^{R}_{(2m-1)\tau_2/2}$,
$i_{(2m-1)\tau_2/2}$ and $i^{r}_{(2m-1)\tau_2/2}$ denote the indices
of $q^R_{{(2m-1)\tau_2/2},\kappa}$, $q_{{(2m-1)\tau_2/2},\kappa}$ and
$q^r_{{(2m-1)\tau_2/2},\kappa}$ respectively. By Lemma 3.1, Lemma 4.1 and the condition in (4.3), we
obtain\begin{eqnarray}i^{R}_{(2m-1)\tau_2/2}&=&i^{E}_{P}(\frac{2}{R^2}I_{2n},2m-1)\nn\\&=&
2\kappa \left(E(\frac{\frac{2}{R^2}(2m-1)\tau_2/2}{2\pi})-1\right)+2(n-\kappa)\left(E(\frac{\frac{2}{R^2}(2m-1)\tau_2/2+\pi}{2\pi})-1\right)
\nn\\&\geq&2nl,\\i(u^{2m-1})&=& i_{(2m-1)\tau_2/2}\geq i^{R}_{(2m-1)\tau_2/2}
\geq 2nl.\end{eqnarray} Hence (4.3) holds. Denote by
the orthogonal splitting $L^2_{\kappa}(0, {(2m-1)\tau_2/2})=E_{-}\oplus
E_{0}\oplus E_{+}$ of $L^2_{\kappa}(0, {(2m-1)\tau_2/2})$ into negative, zero and positive subspaces. Then
we have the following observation: If $V$ is a subspace of
$L^2_{\kappa}(0, {(2m-1)\tau_2/2})$ such
that $q_{{(2m-1)\tau_2/2},\kappa}$ is negative semi-definite, i.e., $v \in
V$ implies $q_{{(2m-1)\tau_2/2},\kappa}(v, v)\leq 0$, then $dimV \leq
dimE_- + dimE_0$. In fact, this is a simple fact of linear algebra:
Let\[pr_-:L^2_{\kappa}(0, {(2m-1)\tau_2/2})=E_{-}\oplus E_{0}\oplus E_{+}\rightarrow E_-\]be the
orthogonal projection. Consider $pr_-|V : V \rightarrow E_-.$ Then $
v \in ker (pr_-|V) $ must belong to $E_{0}$. That is, since
$q_{{(2m-1)\tau_2/2},\kappa}(v, v)> 0, 0\neq v \in E_+$. From\[dimV =
dim~Im (pr_-|V) + dim~ker (pr_-|V)\]we prove our claim.

Let $\epsilon >0$ be small enough such that
$\frac{(2m-1)\tau_2}{2}<(l-\frac{1}{2})\pi (r-\epsilon)^{2}$. If V
is a subspace of $L^2_{\kappa}(0, {(2m-1)\tau_2/2})$ such that $q_{{(2m-1)\tau_2/2},\kappa}\mid V \leq 0$, then
$q^{r-\epsilon}_{{(2m-1)\tau_2/2},\kappa}(v, v)< 0, 0\neq v \in V$. Thus
we have $dimV \leq i^{r-\epsilon}_{(2m-1)\tau_2/2}$. In particular,
by Lemma 3.1, we have
\begin{eqnarray}i(u^{2m-1})+\nu(u^{2m-1})=i_{(2m-1)\tau_2/2}+\nu\left(q_{{(2m-1)\tau_2/2},\kappa}\right)-1\leq
i^{r-\epsilon}_{(2m-1)\tau_2/2}-1.\end{eqnarray}On the other hand,
similarly to (4.5), by the condition in (4.4), we have
\begin{eqnarray}
i^{r-\epsilon}_{(2m-1)\tau_2/2}&=&
2\kappa \left(E(\frac{\frac{2}{(r-\epsilon)^2}(2m-1)\tau_2/2}{2\pi})-1\right)+
2(n-\kappa)\left(E(\frac{\frac{2}{(r-\epsilon)^2}(2m-1)\tau_2/2+\pi}{2\pi})-1\right)\nn\\
&\leq& 2n(l-1).~~~~~\end{eqnarray} Hence,
$i(u^{2m-1})+\nu(u^{2m-1})\leq 2n(l-1)-1$. The proof is complete.

By Theorem 1.1 of \cite{LiZ1}, we have:

{\bf Lemma 4.3.} {\it Assume $\Sigma\in\mathcal {H}_\kappa(2n)$ and $0<r \leq |x|\leq R,~\forall~x \in \Sigma$
with $\frac{R}{r}<\sqrt{2}$. Then there exist at least $n-\kappa$ geometrically distinct P-symmetric closed
characteristics $(\tau_i, y_i)$ on $\Sigma$, where $\tau_i$ is the minimal period of $y_i$,
and the actions $A(\tau_i, y_i)$ satisfy:  \bea \pi r^2\leq A(\tau_i, y_i)\leq \pi R^2, \forall 1\leq i\leq n-\kappa.\eea}

By the proof of the above Lemma and Proposition 2.2, we have

{\bf Theorem 4.4.} {\it Let $\{(\tau_{1}, y_{1}), \cdots ,
(\tau_{n-\kappa}, y_{n-\kappa})\}$ be the P-symmetric closed characteristics found
in the above Lemma. Then we
have\begin{eqnarray}\Psi_{a}^{\prime}(u_{i})=0,~~~ i(u_{i})\leq
2(i-1) \leq i(u_{i})+\nu(u_{i})-1.\end{eqnarray}for $1\leq i \leq
n-\kappa$, where $u_{i}$ is the unique critical point of $\Psi_{a}$
corresponding to $(\tau_{i}, y_{i})$. }

Now we give the proof of the main theorem.

{\bf Proof of Theorem 1.1.} Let $\Sigma\in\mathcal {H}_\kappa(2n)$ which is $(r,R)$-pinched with
$\frac{R}{r}<\sqrt{\frac{5}{3}}$, then by (4.2) we
have\begin{eqnarray}r \leq |x| \leq R,~\forall~x \in
\Sigma\end{eqnarray}Thus by Theorem 4.4, we obtain $n-\kappa$ geometrically
distinct prime P-symmetric closed characteristics $\{(\tau_{1},
y_{1}), \cdots , (\tau_{n-\kappa}, y_{n-\kappa})\}$ such that (4.10) hold.

In the following, we prove $e(\tau_i, y_i)\geq 2n-4\kappa$ for $i=1, n-\kappa$.

Note that we always have $e(\tau_1, y_1)\geq 2n-4\kappa$ by the proof of
Theorem 1 of \cite{DoL2}, $\sqrt{\frac{5}{3}}$-pinching is not necessary.

Note that we can prove $e(\tau_1, y_1)\geq 2n-2\kappa$ under $\sqrt{\frac{3}{2}}$-pinching condition.
In fact, by (4.10), we have $i(u_1)=0$. On the other hand, from (4.9),
we have \bea\frac{(4-1)}{2}A(\tau_1,y_1)\geq\frac{3}{2}\pi
r^{2}> \pi R^{2},\nn\eea where we used the pinching condition
$\frac{R}{r}<\sqrt{\frac{3}{2}}$. By Proposition 4.2, we obtain
\bea i(u_1^3)\geq 2n,\nn\eea
then \bea i(u_1^3)-3i(u_1)\geq 2n.\nn\eea
By Theorem 3.7 (i), we get $e(\tau_1, y_1)\geq 2n-2\kappa$.

Now we prove $e(\tau_{n-\kappa}, y_{n-\kappa})\geq 2n-4\kappa$. In fact, by (4.10), we have
\bea 2(n-\kappa)-1 \leq
i(u_{n-\kappa})+\nu(u_{n-\kappa}).\eea
On the other hand, from (4.9),
we have \bea\frac{(4-1)}{2}A(\tau_{n-\kappa},y_{n-\kappa}) \leq \frac{3}{2}\pi
R^{2}<\frac{5}{2}\pi r^{2}=(3-\frac{1}{2})\pi r^2.\nn\eea
By Proposition 4.2, we obtain
\bea i(u_{n-\kappa}^3)+\nu(u_{n-\kappa}^3)\leq 4n-1,\eea
Combining (4.12) with (4.13), we obtain \bea i(u_{n-\kappa}^3)+\nu(u_{n-\kappa}^3)-3(i(u_{n-\kappa})+\nu(u_{n-\kappa}))
\leq 6\kappa+2-2n, \eea
then $e(\tau_{n-\kappa}, y_{n-\kappa})\geq 2n-4\kappa$ follows from Theorem 3.7 (ii). The proof is complete.

{\bf Acknowledgements.} The authors would like to sincerely thank the anonymous referee for his/her
careful reading of the manuscript and valuable comments.

\bibliographystyle{abbrv}

\begin{thebibliography}{*}
\bibitem[DoL1]{DoL1} Y. Dong, Y. Long, Closed characteristics on partially symmetric compact convex
hypersurfaces in ${\bf R}^{2n}$. {\it J. Differential Equations.} 196
(2004), 226-248.
\bibitem[DoL2]{DoL2} Y. Dong, Y. Long, Stable closed characteristics on partially symmetric compact convex
hypersurfaces. {\it J. Differential Equations.} 206 (2004), 265-279.
\bibitem[Eke1]{Eke1} I. Ekeland, Convexity Methods in Hamiltonian Mechanics. Springer-Verlag. Berlin. 1990.
\bibitem[EkH1]{EkH1} I. Ekeland and H. Hofer, Convex Hamiltonian energy
surfaces and their closed trajectories, {\it Comm. Math. Phys.} 113
(1987), 419-467.
\bibitem[EkL1]{EkL1} I. Ekeland and J. Lasry, On the number of periodic trajectories for a Hamiltonian flow on a
convex energy surface. {\it Ann. of Math.} 112 (1980), 283-319.
\bibitem[EkL2]{EkL2} I. Ekeland and L. Lassoued, Multiplicit\'e des
trajectoires ferm\'ees d'un syst\'eme hamiltonien sur une
hypersurface d'energie convexe. {\it Ann. IHP. Anal. non
Lin\'eaire.} 4 (1987), 1-29.
\bibitem[FaR1]{FaR1} E. Fadell and P. Rabinowitz, Generalized
cohomological index theories for Lie group actions with an application
to bifurcation equations for Hamiltonian systems. {\it Invent. Math.
} 45 (1978), 139-174.
\bibitem[Gir1]{Gir1} M. Girardi, Multiple orbits for Hamiltonian
systems on starshaped ernergy surfaces with symmetry. {\it Ann. IHP.
Analyse non lin\'eaire.} 1 (1984), 285-294.
\bibitem[HuS1]{HuS1} X. Hu and S. Sun, Index and stability of symmetric periodic orbits in Hamiltoinan systems
with application to figure-eight orbits, {\it Commun. Math. Phys.} 290 (2009), 737-777.
\bibitem[Liu1]{Liu1} H. Liu, Stability of symmetric closed characteristics on symmetric compact convex hypersurfaces
in ${\bf R}^{2n}$ under a pinching condition. {\it Acta Mathematica
Sinica, English Series,} 28 (2012), 885-900.
\bibitem[Liu2]{Liu2} H. Liu, Multiple P-invariant closed characteristics on partially symmetric compact convex hypersurfaces in ${\bf R}^{2n}$.
{\it Cal. Variations and PDEs.} 49 (2014), 1121-1147.
\bibitem[LiZ1]{LiZ1} H. Liu and D. Zhang, On the number of P-invariant closed characteristics on partially symmetric compact convex
hypersurfaces in ${\bf R}^{2n}$. {\it Science China Mathematics.} 2015, 58, doi: 10.1007/s11425-014-4903-2.
\bibitem[LLZ1]{LLZ1} C. Liu, Y. Long and C. Zhu, Multiplicity of closed characteristics on symmetric convex
hypersurfaces in ${\bf R}^{2n}$. {\it Math. Ann.} 323 (2002),
201-215.
\bibitem[Lon1]{Lon1} Y. Long, Index Theory for Symplectic Paths
with Applications. Progress in Math. 207, Birkh\"auser. Basel. 2002.
\bibitem[LoZ1]{LoZ1} Y. Long and C. Zhu, Closed characteristics on compact convex hypersurfaces in ${\bf R}^{2n}$,
{\it Ann. of Math.} 155 (2002), 317-368.
\bibitem[Rab1]{Rab1} P. Rabinowitz, Periodic solutions of Hamiltonian systems, {\it Comm. Pure. Appl.
Math.} 31 (1978), 157-184.
\bibitem[Szu1]{Szu1} A. Szulkin, Morse theory and existence of periodic solutions
of convex Hamiltonian systems, {\it Bull. Soc. Math. France.} 116
(1988), 171-197.
\bibitem[Wan1]{Wan1} W. Wang, Closed trajectories on symmetric convex
Hamiltonian energy surfaces. {\it Discrete Contin. Dyn. Syst.} 32 (2012), no. 2, 679-701.
\bibitem[Wei1]{Wei1} A. Weinstein, Periodic orbits for convex Hamiltonian systems, {\it Ann. of Math.} 108
(1978), 507-518.
\bibitem[WHL1]{WHL1} W. Wang, X. Hu and Y. Long, Resonance identity, stability and multiplicity of closed
characteristics on compact convex hypersurfaces, {\it Duke Math. J.}
139 (2007), 411-462.
\bibitem[Zha1]{Zha1} D. Zhang, P-cyclic symmetric closed characteristics on compact convex P-cyclic symmetric hypersurface in ${\bf R}^{2n}$.
{\it Discrete Contin. Dyn. Syst.} 33 (2013), no. 2, 947-964.

\end{thebibliography}

\end{document}